\documentclass[12pt]{amsart}
\usepackage{amsmath,amssymb}
\usepackage{fullpage}

\theoremstyle{plain}
\newtheorem{thm}{Theorem}[section]

\newtheorem{lem}[thm]{Lemma}
\theoremstyle{definition}
\newtheorem{defn}[thm]{Definition}
\newtheorem{exa}[thm]{Example}
\newtheorem{rem}[thm]{Remark}

\numberwithin{equation}{section}

\begin{document}
\title[Nonlocal Caputo FBVPs]{Nontrivial solutions of systems of nonlocal Caputo fractional BVPs}




\author[G. Infante]{Gennaro Infante}
\address{Gennaro Infante, Dipartimento di Matematica e Informatica, Universit\`{a} della
Calabria, 87036 Arcavacata di Rende, Cosenza, Italy}%
\email{gennaro.infante@unical.it}%

\author[S. Rihani]{Samira Rihani}%
\address{Samira Rihani, Faculty of Mathematics, USTHB, El Alia Bab Ezzouar
Algiers 16111, Algeria}%
\email{maths\_samdz@yahoo.fr}%

\begin{abstract}
We discuss the existence, non-existence and multiplicity of nontrivial solutions for systems of Caputo fractional differential equations subject to nonlocal boundary conditions. 
Our methodology relies on classical fixed point index and we make use of recent results by Infante and Pietramala.
\end{abstract}

\maketitle

\section{Introduction}

Nieto and Pimentel~\cite{nieto-pim}, motivated by earlier work of Infante and Webb~\cite{gijwnodea} on nonlocal problems for ODEs,  studied the existence of at least one positive solution of the fractional differential equation
\begin{equation}\label{niepim-eq}
{}^C\!D^{\alpha} u(t)+f(t, u(t))=0, \ t\in(0, 1),
\end{equation}
subject to the nonlocal boundary conditions (BCs)
\begin{equation}\label{niepim-bcs}
u'(0)=0,\ \beta{}^C\!D^{\alpha-1}u(1)+u(\eta)=0,
\end{equation}
where $1<\alpha\leq 2$, ${}^C\!D^{\alpha}$ denotes the Caputo fractional derivative of order $\alpha$, $\beta>0$, $0\leq \eta\leq 1$ and $f$ is continuous. The study of Nieto and Pimentel was continued by Cabada and Infante~\cite{ac-gi} who, by means of fixed point index theory, studied the existence of \emph{multiple} positive solutions of~\eqref{niepim-eq} under some BCs that involve Riemann-Stieltjes integrals and cover the ones in~\eqref{niepim-bcs} as a special case. 

Note that the BCs that occur in~\eqref{niepim-bcs} are of \emph{nonlocal} type. The study of nonlocal BCs goes back, as far as we know, to Picone~\cite{Picone} and has been widely developed during the years. We refer the reader to the reviews by
Whyburn~\cite{Whyburn}, Conti~\cite{Conti}, Ma~\cite{rma},  Ntouyas~\cite{sotiris} and \v{S}tikonas~\cite{Stik} and the papers by Karakostas and Tsamatos~\cite{kttmna, ktejde} and by Webb and Infante~\cite{jwgi-lms}. We mention that Costabile and Napoli~\cite{cos} also contributed to the study of nonlocal problems in the context of ODEs under a numerical point of view.

The problem of existence of solutions for \emph{systems} of Caputo fractional differential equations under a variety of   
BCs has been investigated by a number of authors; for example local BCs have been investigated by Khan and  ur~Rehman~\cite{khan-re} and Lan and Lin~\cite{lanfrac} and nonlocal BCs have been studied by Ahmad and Nieto~\cite{ahma-nie}, ur~Rehman and co-authors~\cite{re-kh-el} and Zhao and Gong~\cite{zha-go}.

In this paper we discuss the existence, non-existence and multiplicity of \emph{nontrivial} solutions for the system of fractional differential equations
\begin{gather}
\begin{aligned}\label{1syst}
{}^C\!D^{\alpha_1}u(t) + f_1(t,u(t),v(t)) =& 0, \ t \in (0,1), \\
{}^C\!D^{\alpha_2}v(t) + f_2(t,u(t),v(t)) =& 0, \ t \in (0,1),%
\end{aligned}
\end{gather}
subject the BCs
\begin{gather}
\begin{aligned}\label{1BC}
u'(0)=0,\ \beta_1{}^C\!D^{\alpha_1-1}u(1)+u(\eta_1)=0, \\
v'(0)=0,\ \beta_2{}^C\!D^{\alpha_2-1}v(1)+v(\eta_2)=0,
\end{aligned}
\end{gather}
where $1<\alpha_1, \alpha_2\leq 2$, $\beta>0$, $0\leq \eta_1, \eta_2 \leq 1$ and $f_1, f_2$ are continuous.
In order to do this, we associate to the system \eqref{1syst}-\eqref{1BC} a
system of Hammerstein integral equations of the type
\begin{gather}
\begin{aligned}\label{syst}
u(t)=\int_{0}^{1}k_1(t,s)f_1(s,u(s),v(s))\,ds, \\
v(t)= \int_{0}^{1}k_2(t,s)f_2(s,u(s),v(s))\,ds,%
\end{aligned}
\end{gather}
and we make use of recent results of Infante and Pietramala~\cite{gipp}, that rely 
on the fixed point index theory.
We work in a suitable cone in a product space of continuous functions, where functions are allowed to \emph{change sign}. 

Our results complement those in~\cite{nieto-pim}, providing multiplicity and nonexistence results in the context of nontrivial solutions.

\section{Existence and nonexistence results}\label{secham}
We begin by recalling some results from~\cite{gipp} and we make the following assumptions on the terms that occur in the system~\eqref{syst}.
\begin{itemize}
\item For every $i=1,2$, $f_i: [0,1]\times (-\infty,\infty)\times (-\infty,\infty) \to
[0,\infty)$ is continuous.
{}
\item For every $i=1,2$, $k_i:[0,1]\times [0,1]\to (-\infty,\infty)$  is continuous.
{}
\item For every $i=1,2$, there exist a subinterval $[a_i,b_i] \subseteq
[0,1] $, a function $\Phi_i \in L^{\infty}[0,1]$, and a constant $c_{i} \in
(0,1]$, such that
\begin{align*}
|k_i(t,s)|\leq \Phi_i(s) \text{ for } &t \in [0,1] \text{ and a.\,e.}\, s\in [0,1], \\
k_i(t,s) \geq c_{i}\Phi_i(s) \text{ for } &t\in [a_i,b_i] \text{ and a.\,e.} \, s \in [0,1].
\end{align*}%
{}
\item For every $i=1,2$, we have $\int_{a_i}^{b_i} \Phi_i(s)\,ds >0$.
\end{itemize}

We work in the space $C[0,1]\times C[0,1]$ endowed with the norm
\begin{equation*}
\| (u,v)\| :=\max \{\| u\| _{\infty },\| v\| _{\infty }\},
\end{equation*}%
where $\| w\| _{\infty}:=\max \{| w(t)|,t\in [0,1] \}$, define the sets
\begin{equation*}
\tilde{K_{i}}:=\{w\in C[0,1]:\min_{t\in \lbrack a_{i},b_{i}]}w(t)\geq c_{i}\|
w\| _{\infty }\},
\end{equation*}%
 and consider the cone $K$
in $C[0,1]\times C[0,1]$ defined by
\begin{equation*}
\begin{array}{c}
K:=\{(u,v)\in \tilde{K_{1}}\times \tilde{K_{2}}\}.%
\end{array}%
\end{equation*}

By a \emph{nontrivial} solution of the system \eqref{syst} we mean a solution
$(u,v)\in K$ of \eqref{syst} such that $\|(u,v)\|\neq 0$.
Note that the functions in $\tilde{K_{i}}$ are positive on the
sub-interval $[a_{i},b_{i}]$ but are allowed to change sign in $[0,1]$. This type of cone has been introduced by  Infante and Webb in \cite{gijwjiea}
and is similar to
a cone of \emph{non-negative} functions
firstly used by Krasnosel'ski\u\i{}, see e.g. \cite{krzab}, and D.~Guo, see e.g. \cite{guolak}. 

Note that, under the assumptions above (see Lemma~2.1 of ~\cite{gipp}), the integral operator
\begin{gather}
\begin{aligned} \label{opT}
T(u,v)(t):=& 
\left(
\begin{array}{c}
 \int_{0}^{1}k_1(t,s)f_1(s,u(s),v(s))\,ds\\
\int_{0}^{1}k_2(t,s)f_2(s,u(s),v(s))\,ds%
\end{array}
\right)
\end{aligned}
\end{gather}
leaves the cone $K$ invariant and is compact.

We use the following open bounded sets (relative to $K$):
\begin{equation*}
K_{\rho_1,\rho_2} = \{ (u,v) \in K : \|u\|_{\infty}< \rho_1\ \text{and}\ \|v\|_{\infty}< \rho_2\},
\end{equation*}
and
\begin{equation*}
V_{\rho_1,\rho_2} =\{(u,v) \in K: \min_{t\in [a_1,b_1]}u(t)<\rho_1\ \text{and}\ \min_{t\in
[a_2,b_2]}v(t)<\rho_2\}.
\end{equation*}
Note that
$K_{\rho_1,\rho_2}\subset V_{\rho_1,\rho_2}\subset K_{\rho_1/c_1,\rho_2/c_2}$; this is a key property used in order to prove the multiplicity results. 
When $\rho_1=\rho_2=\rho$ we write $K_{\rho}$ and $V_{\rho}$.
The set $V_\rho$ (in the context of systems) was introduced by Infante and Pietramala in~\cite{gipp-ns} and is equal to the set
called $\Omega^{\rho /c}$ in~\cite{df-gi-do}. $\Omega^{\rho /c}$ is an extension to the case of systems of a set given by Lan~\cite{lan}. 

The next Lemma summarizes some sufficient conditions from~\cite{gipp} regarding the index computations.
\begin{lem}\label{index} 
The  following hold.
\begin{enumerate}
\item Assume that 
\begin{enumerate}
\item[$(\mathrm{I}_{\rho_1,\rho_2 }^{1})$] \label{EqB} there exist $\rho_1,\rho_2 >0$
such that for every $i=1,2$
\begin{equation}\label{eqmestt}
 f_i^{\rho_1,\rho_2}  < m_i
\end{equation}{}
where
\begin{equation*}
f_{i}^{\rho_1,\rho_2}=\sup \Bigl\{\frac{f_{i}(t,u,v)}{\rho_i }:\;(t,u,v)\in
\lbrack 0,1]\times [ -\rho_1,\rho_1 ]\times [ -\rho_2,\rho_2 ]\Bigr\} \end{equation*}
and
\begin{equation*}
 \frac{1}{m_{i}}=\sup_{t\in \lbrack 0,1]}\int_{0}^{1}\vert k_{i}(t,s)\vert \,ds.
\end{equation*}{}
\end{enumerate}{}
Then the fixed point index of $T$ relative to $K_{\rho}$, $i_{K}(T,K_{\rho_1,\rho_2})$ is equal to $1$.
\item Assume that 
\begin{enumerate}
\item[$(\mathrm{I}^{0}_{\rho_1,\rho_2})$]  there exist $\rho_1,\rho_2>0$ such that for every $i=1,2$
\begin{equation}\label{eqMest}
f_{i,(\rho_1, \rho_2)} > M_i,
\end{equation}{}
where
\begin{align*}
 f_{1,({\rho_1,\rho_2 })}=& \inf \Bigl\{ \frac{f_1(t,u,v)}{ \rho_1}:\; (t,u,v)\in [a_1,b_1]\times[\rho_1,\rho_1/c_1]\times[-\rho_2/c_2, \rho_2/c_2]\Bigr\},\\
 f_{2,({\rho_1,\rho_2 })}=& \inf \Bigl\{ \frac{f_2(t,u,v)}{ \rho_2}:\; (t,u,v)\in [a_2,b_2]\times[-\rho_1/c_1,\rho_1/c_1]\times[\rho_2, \rho_2/c_2]\Bigr\},\\
    \frac{1}{M_i}=& \inf_{t\in
[a_i,b_i]}\int_{a_i}^{b_i} k_i(t,s) \,ds.
\end{align*}
\end{enumerate}
Then $i_{K}(T,V_{\rho_1,\rho_2})=0$.
\item Assume that
\begin{enumerate}
\item[$(\mathrm{I}^{0}_{\rho_1,\rho_2})^{\star}$] there exist $\rho_1,\rho_2>0$ such that for some $i\in\{1,2\}$ we have
\begin{equation}\label{eqMest1}
f^*_{i,(\rho_1, \rho_2)}>M_i,
\end{equation}{}
\end{enumerate}
where
\begin{equation*}
f^*_{1,(\rho_1,{\rho_2})}=\inf \Bigl\{ \frac{f_1(t,u,v)}{ \rho_1}:\; (t,u,v)\in [a_1,b_1]\times[0,\rho_1/c_1]\times[-\rho_2/c_2, \rho_2/c_2]\Bigr\}.
\end{equation*}
\begin{equation*}
f^*_{2,(\rho_1,{\rho_2})}=\inf \Bigl\{ \frac{f_2(t,u,v)}{ \rho_2}:\; (t,u,v)\in [a_2,b_2]\times[-\rho_1/c_1,\rho_1/c_1]\times[0, \rho_2/c_2]\Bigr\}.
\end{equation*}
Then $i_{K}(T,V_{\rho_1,\rho_2})=0$.
\end{enumerate}
\end{lem}
\begin{rem}
In the previous Lemma the condition $(\mathrm{I}^{0}_{\rho_1,\rho_2})^{\star}$ allows the nonlinearities to have a different growth, at the cost of having to deal with a larger domain. Nonlinearities with different growths were considered, with different approaches, in~\cite{chzh1, precup1, precup2, ya1}.
\end{rem}

Using the Lemma~\ref{index} it is possible to prove, by means of fixed point index, a result regarding the existence of at least one, two or three nontrivial solutions. Note that it is also possible to state a results for four or more nontrivial solutions, in the line of the paper~\cite{kljdeds}.

\begin{thm}\label{mult-sys}\cite{gipp}
The system \eqref{syst} has at least one nontrivial solution
in $K$ if one of the following conditions holds.

\begin{enumerate}

\item[$(S_{1})$]  For $i=1,2$ there exist $\rho _{i},r _{i}\in (0,\infty )$ with $\rho
_{i}/c_i<r _{i}$ such that $(\mathrm{I}_{\rho _{1},\rho_2}^{0})\;\;[\text{or}\;(%
\mathrm{I}_{\rho _{1},\rho_2}^{0})^{\star }]$, $(\mathrm{I}_{r _{1},r_2}^{1})$ hold.

\item[$(S_{2})$] For $i=1,2$ there exist $\rho _{i},r _{i}\in (0,\infty )$ with $\rho
_{i}<r _{i}$ such that $(\mathrm{I}_{\rho _{1},\rho_2}^{1}),\;\;(\mathrm{I}%
_{r _{1},r_2}^{0})$ hold.
\end{enumerate}

The system \eqref{syst} has at least two nontrivial solutions in $K$ if one of the following conditions holds.

\begin{enumerate}

\item[$(S_{3})$] For $i=1,2$ there exist $\rho _{i},r _{i},s_i\in (0,\infty )$
with $\rho _{i}/c_i<r_i <s _{i}$ such that $(\mathrm{I}_{\rho
_{1},\rho_2}^{0})$, $[\text{or}\;(\mathrm{I}_{\rho _{1},\rho_2}^{0})^{\star }],\;\;(%
\mathrm{I}_{r _{1},r_2}^{1})$ $\text{and}\;\;(\mathrm{I}_{s _{1},s_2}^{0})$
hold.

\item[$(S_{4})$] For $i=1,2$ there exist $\rho _{i},r _{i},s_i\in (0,\infty )$
with $\rho _{i}<r _{i}$ and $r _{i}/c_i<s _{i}$ such that $(\mathrm{I}%
_{\rho _{1},\rho_2}^{1}),\;\;(\mathrm{I}_{r _{1},r_2}^{0})$ $\text{and}\;\;(\mathrm{I%
}_{s _{1},s_2}^{1})$ hold.
\end{enumerate}

The system \eqref{syst} has at least three nontrivial solutions in $K$ if one
of the following conditions holds.

\begin{enumerate}
\item[$(S_{5})$] For $i=1,2$ there exist $\rho _{i},r _{i},s_i,\sigma_i\in
(0,\infty )$ with $\rho _{i}/c_i<r _{i}<s _{i}$ and $s _{i}/c_i<\sigma
_{i}$ such that $(\mathrm{I}_{\rho _{1},\rho_2}^{0})\;\;[\text{or}\;(\mathrm{I}%
_{\rho _{1},\rho_2}^{0})^{\star }],$ $(\mathrm{I}_{r _{1},r_2}^{1}),\;\;(\mathrm{I}%
_{s_1,s_2}^{0})\;\;\text{and}\;\;(\mathrm{I}_{\sigma _{1},\sigma_2}^{1})$ hold.

\item[$(S_{6})$] For $i=1,2$ there exist $\rho _{i},r _{i},s_i,\sigma_i\in
(0,\infty )$ with $\rho _{i}<r _{i}$ and $r _{i}/c_i<s _{i}<\sigma _{i}$
such that $(\mathrm{I}_{\rho _{1},\rho_2}^{1}),\;\;(\mathrm{I}_{r
_{1},r_2}^{0}),\;\;(\mathrm{I}_{s _{1},s_2}^{1})$ $\text{and}\;\;(\mathrm{I}%
_{\sigma _{1},\sigma_2}^{0})$ hold.
\end{enumerate}
\end{thm}

The following result provides some sufficient conditions for the non-existence of nonzero solutions.
\begin{thm}\cite{gipp}\label{non-exist}
Assume that one of the following conditions holds: 
\begin{enumerate}
\item For  $i=1,2$, 
\begin{equation}\label{cond1}
f_i(t,u_1,u_2)<m_i|u_i| \text{         for every        } t\in [0,1] \text{        and        } u_i\neq 0.
\end{equation}
\item For  $i=1,2$, 
\begin{equation}\label{cond2}
f_i(t,u_1,u_2)>M_i u_i \text{         for every        } t\in [a_i,b_i] \text{        and        }  u_i>0.
\end{equation}
\item There exists $i=1,2$ such that \eqref{cond1} is verified for $f_i$ and for  $j\neq i$  condition \eqref{cond2}  is verified for $f_j$.
\end{enumerate}
Then there is no non-trivial solution of the system \eqref{syst} in $K$. 
\end{thm}

We now turn our attention back to the fractional system \eqref{1syst}-\eqref{1BC} and illustrate the applicability of the above results. 
Firstly we recall the definition of the Caputo derivative, for its properties we refer to the books~\cite{anast, dieth, pod, samk}. 

\begin{defn}
For a function $y:[0, +\infty)\to \mathbb{R}$,
the  Caputo  derivative of fractional order
  $\alpha>0$  is given by
$$
{}^C\!D^\alpha
y(t)=\frac{1}{\Gamma(n-\alpha)}\int^t_0\frac{y^{(n)}(s)}{(t-s)^{\alpha+1-n}}\,ds, \quad n=[\alpha] + 1,
$$
where $\Gamma$ denotes the Gamma function, that is
$$\Gamma (s)=\int_{0}^{+\infty}x^{s-1}e^{-x}dx,$$
and $[\alpha]$ denotes the integer part of a number $\alpha$.
\end{defn}
Observe that, by direct calculation (see Lemma 2.4 of \cite{nieto-pim}), the solution of the linear equation
$$
{}^C\!D^{\alpha} w(t)+y(t)=0,
$$
under the BCs 
$$
w'(0)=0,\ \beta{}^C\!D^{\alpha-1}w(1)+w(\eta)=0,
$$
can be written in the form
$$
w(t)=\beta \int_{0}^{1}y(s)ds
+\int_{0}^{\eta}\frac{(\eta-s)^{\alpha-1}}{\Gamma(\alpha)}y(s)ds -\int_{0}^{t}
\frac{(t-s)^{\alpha-1}}{\Gamma(\alpha)}y(s)ds. 
$$

Therefore the solution of the system \eqref{1syst}-\eqref{1BC} is given by
\begin{gather}
\begin{aligned}\label{systapp}
u(t)=\int_{0}^{1}k_1(t,s)f_1(s,u(s),v(s))\,ds, \\
v(t)= \int_{0}^{1}k_2(t,s)f_2(s,u(s),v(s))\,ds,%
\end{aligned}
\end{gather}
where 
$$
k_i(t,s)=\beta_i +\frac{1}{\Gamma(\alpha_i)}
\begin{cases} (\eta_i-s)^{\alpha_i-1},  s\leq \eta_i,\\ 
0,\hfill s>\eta_i,
\end{cases}
-\frac{1}{\Gamma(\alpha_i)}
\begin{cases}(t-s)^{\alpha_i-1},  s\leq t,\\ 
0,\hfill  s>t.
\end{cases}
$$
Here we focus on the case $$\beta_i  \Gamma(\alpha_i) < (1-\eta_i )^{\alpha_i-1},$$
and seek solutions that are positive on a subinterval of $[0,1]$ and are allowed to change sign elsewhere.
In this case the interval
$[a_i, b_i]$ can be chosen equal to $[0,b_i]$, where $b_i$ is such that $\eta_i\leq b_i<1$ and  
$\beta_i  \Gamma(\alpha_i)>(b_i-\eta_i )^{\alpha_i-1}$.\medskip

Upper and lower bounds for $k_i(t,s)$ were given in \cite{nieto-pim}, as follows:
\begin{eqnarray}
\Phi_i(s)=\begin {cases}
\frac{(1-\eta_i)^{\alpha_i-1}}{\Gamma(\alpha_i)}-\beta_i, & s>\eta_i, \\
\beta_i+\frac{(\eta_i-s)^{\alpha_i-1}}{\Gamma(\alpha_i)}, & s\leq \eta_i,
	\end{cases}
\end{eqnarray}
$$
c_i=\min \Bigl\{\frac{\beta_i\Gamma(\alpha_i)-(b_i-\eta_i)^{\alpha_i-1}}{(1-\eta_i)^{\alpha_i-1}-\beta_i\Gamma(\alpha_i)}
,\, \frac{\beta_i\Gamma(\alpha_i)-(b_i-\eta_i)^{\alpha_i-1}}{\beta_i\Gamma(\alpha_i)+\eta_i^{\alpha_i-1}}
\Bigr\}. 
$$

Thus we work in the cone \begin{equation*}
\begin{array}{c}
K:=\{(u,v)\in \tilde{K_{1}}\times \tilde{K_{2}}\},%
\end{array}%
\end{equation*}
where
\begin{equation*}
\tilde{K_{i}}:=\{w\in C[0,1]:\min_{t\in \lbrack 0,b_{i}]}w(t)\geq c_{i}\|
w\| _{\infty }\}.
\end{equation*}

Theorems \ref{mult-sys} and \ref{non-exist} can be applied to the system \eqref{1syst}-\eqref{1BC}, provided that the nonlinearities have a suitable growth. We stress that the constants that occur in our theory can either be computed or estimated. In the next example we show precisely this situation.
The numbers are rounded to the third decimal place, unless exact.
\begin{exa}
Consider the data
$$\alpha_1= 3/2,\ \beta_1= 1/5,\ \eta_1 = 3/4,$$
$$\alpha_2= 5/4,\ \beta_2 = 2/5,\  \eta_2 = 2/3. $$
In this case, one may use the intevals
$$[a_1,b_1]=[0, 31/40],\ [a_2,b_2]=[0, 41/60],$$
and this choice leads to 
$$c_1=0.018,\ c_2=0.002.$$
We make use of quantities
$$
\frac{1}{\hat{M}_{i}}:=\int_{a_i}^{b_i} c\Phi_i(s) \,ds,\ \frac{1}{\hat{m}_i}:=\int_{0}^{1} \Phi_i(s) \,ds
$$
and note that 
$$M_i \leq \hat{M}_i \ \text{and}\ m_i \geq \hat{m}_i.$$
With the data above we obtain
$$ \hat{M}_1=84.192\quad \text{and}\quad \hat{m}_1=1.370.$$
$$ \hat{M}_2=482.545\quad \text{and}\quad \hat{m}_2=1.058.$$
\end{exa}
\section*{Acknowledgments}
G. Infante was partially supported by G.N.A.M.P.A. - INdAM (Italy).

\end{document}